\documentclass{article}
\usepackage{amssymb,amsmath,amsthm,graphicx}

\newtheorem{theorem}{Theorem}

\theoremstyle{definition}
\newtheorem{definition}{Definition}
\newtheorem{example}{Example}
\newtheorem{remark}{Remark}

\date{}

\title{\Large \textbf{Finite Type Enhancements}}

\author{Sam Nelson\footnote{Email: knots@esotericka.org. 
Partially supported by Simons Foundation Collaboration Award 316709}}

\begin{document}
\maketitle

\begin{abstract}
We enhance the biquandle counting invariant using elements of truncated 
biquandle-labeled Polyak algebras. These finite type enhancements reduce 
to the finite type enhancements defined in \cite{GPV} for the trivial 
biquandle of one element and determine (but are not determined by) the 
biquandle counting invariant for general biquandles. Unlike the unlabeled
case, biquandle labeled finite type invariants of degree 1 are nontrivial
and are related to biquandle cocycle invariants.
\end{abstract}

\parbox{5.5in}{\textsc{Keywords:} Finite type invariants, Virtual knots, 
Biquandles, Enhancements of counting invariants

\smallskip

}

\section{\large\textbf{Introduction}}\label{I}

\textit{Biquandles} are a type of algebraic structure with axioms motivated 
by knot theory. First introduced in \cite{FRS1}, biquandles have been further
developed in many recent works such as \cite{FJK,FRS,KR}. Given a finite 
biquandle $X$, the set of $X$-labelings of a tame oriented knot or link diagram
is a finite set whose elements can be interpreted as biquandle homomorphisms
from the \textit{fundamental biquandle} of the knot or link to $X$. The 
fundamental biquandle of a knot determines the fundamental quandle of the knot,
a complete invariant for classical oriented knots up to reflection, and is 
conjectured to be a complete invariant of virtual knots up to reflection 
\cite{FJK}. 

Any invariant $\phi$  of $X$-labeled diagrams determines an invariant of 
oriented knots and links, namely the multiset of $\phi$ values over the set
of biquandle labelings of the knot or link. When $\phi$ takes numerical values,
we can convert the multiset into a polynomial for simplicity. Such an invariant
is called an \textit{enhancement} of the biquandle counting invariant. 
Enhancements are tools for extracting information from the fundamental 
biquandle of a knot or link in a practical way.

Signed Gauss codes and Gauss diagrams were described in \cite{K}. In 
\cite{GPV} signed Gauss diagrams were used to define \textit{finite type 
invariants} using an algebra generated by virtual knot diagrams. Unlike the 
unlabeled Polyak algebra, we find that the space of degree 1 biquandle labeled
finite type invariants is generally nontrivial.

In this paper we define finite type invariants for biquandle-labeled knot 
diagrams, obtaining new enhancements of the biquandle counting invariant.
For every pair of finite biquandle $X$ and positive integer $n$, we obtain
a finite dimensional algebra $\mathcal{P}_n^X$, each element of which defines 
an enhancement of the biquandle counting invariant. The paper is organized
as follows. In Section \ref{B} we review the basics of biquandles. In Section
\ref{A} we define biquandle-labeled Polyak algebras and compute some
examples. In Section \ref{E} we define finite type enhancements and demonstrate 
that the enhancements are proper, i.e. not determined by the biquandle 
counting invariants; indeed, we give an example of a degree 1 finite type 
enhancement which can show non-classicality of a virtual knot. In Section
\ref{MPA} we extend our previous work to the case of multicomponent Polyak
algebras and give an example of a finite type enhancement which determines
linking number. We conclude in Section \ref{Q} with some
questions for future research.


This work was partially supported by funding from the Pomona College Department 
of Mathematics, Claremont McKenna College and the Simons Foundation. The author
wishes to thank Pomona College student Selma Paketci for her contributions
to this project.

\section{\large\textbf{Biquandles}}\label{B}

In this section we review the basics of biquandles and the biquandle counting
invariant. We begin with a definition:

\begin{definition}
Let $X$ be a set. A \textit{biquandle structure} on $X$ is a pair of binary
operations denoted by $(x,y)\mapsto x^y, x_y$ such that for all $x,y,z$ we
have
\begin{itemize}
\item[(i)] $x^x=x_x$,
\item[(ii)] The maps $x\to x^y$, $y\to y_x$ and
$(x,y)\mapsto(y_x,x^y)$ are invertible, and
\item[(iii)] The operations satisfy the \textit{exchange laws}
\[\begin{array}{rcl}
(x^y)^{(z^y)} & = & (x^z)^{(y_z)} \\
(x^y)_{(z^y)} & = & (x_z)^{(y_z)} \\
(x_y)_{(z_y)} & = & (x_z)_{(y^z)}.
\end{array}\]
A biquandle in which $y_x=y$ for all $x,y$ is a \textit{quandle}.
\end{itemize}
\end{definition}

\begin{example}
Let $n\in \mathbb{Z}$. A group $G$ is a biquandle (indeed, a quandle) under 
the operations
\[x^y=y^{-n}xy^n\quad\mathrm{and}\quad  y_x=y.\]
Such a quandle is known as the \textit{$n$-fold conjugation quandle}
of the group $G$, denoted $\mathrm{Conj}_n(G)$.
If $n=0$ (so $x_y=x^y=x$ for all $X$), we have a \textit{trivial biquandle}.
\end{example}

\begin{example}
Let $X$ be any set and let $\sigma:X\to X$ be a bijection.
Then $X$ is a biquandle under 
\[x^y=\sigma(x)\quad \mathrm{and}\quad y_x=\sigma(y).\] 
Such a biquandle is known as a \textit{constant action biquandle}.
\end{example}

\begin{example}
Let $X$ be any module over the ring 
$\overline{\Lambda}=\mathbb{Z}[t^{\pm 1},s^{\pm 1}]$. Then $X$ is a biquandle
under the operations
\[x^y=tx+(1-s^{-1}t)y\quad \mathrm{and}\quad y_x=s^{-1}y\]
known as an \textit{Alexander biquandle}. If $s=1$ then $X$ is an 
\textit{Alexander quandle}.
\end{example}

\begin{example}
For a more concrete example, let $X=\mathbb{Z}_3$ and let $t=1$ and $s=2.$
Then $s^{-1}=2$, $x^y=tx+(1-s^{-1}t)y=x+2y$ and $y_x=s^{-1}y=2y$. If we write
$\mathbb{Z}_3=\{0,1,2\}$ then we have operation tables
\[\begin{array}{c|ccc}
x^y & 0 & 1 & 2 \\ \hline
0 & 0 & 2 & 1 \\
1 & 1 & 0 & 2 \\
2 & 2 & 1 & 0
\end{array}
\quad
\begin{array}{c|ccc}
x_y & 0 & 1 & 2 \\ \hline
0 & 0 & 0 & 0 \\
1 & 2 & 2 & 2 \\
2 & 1 & 1 & 1.
\end{array}
\]
\end{example}

More generally, for any finite set $X=\{x_1,\dots, x_n\}$, we can specify a 
biquandle structure on $X$ with a pair of operation tables such that the
axioms are satisfied; for brevity we generally write the operation tables
as an $n\times 2n$ block matrix with entries in $\{1,2,\dots, n\}$.

\begin{example}
The biquandle structure on the set $X=\{x_1,x_2,x_3\}$ with operation tables
\[\begin{array}{c|ccc}
x^y & x_1 & x_2 & x_3 \\ \hline
x_1 & x_1 & x_3 & x_2 \\
x_2 & x_3 & x_2 & x_1 \\
x_3 & x_2 & x_1 & x_3
\end{array}
\quad
\begin{array}{c|ccc}
x_y & x_1 & x_2 & x_3 \\ \hline
x_1 & x_1 & x_1 & x_1 \\
x_2 & x_2 & x_2 & x_2 \\
x_3 & x_3 & x_3 & x_3
\end{array}\]
can be expressed compactly with the biquandle matrix
\[
\left[
\begin{array}{rrr|rrr}
1 & 3 & 2 & 1 & 1 & 1 \\
3 & 2 & 1 & 2 & 2 & 2 \\
2 & 1 & 3 & 3 & 3 & 3
\end{array}\right].
\]
\end{example}

The biquandle axioms are the conditions required to make labelings of the 
semiarcs of oriented knot and link diagrams according to the rules
\[\includegraphics{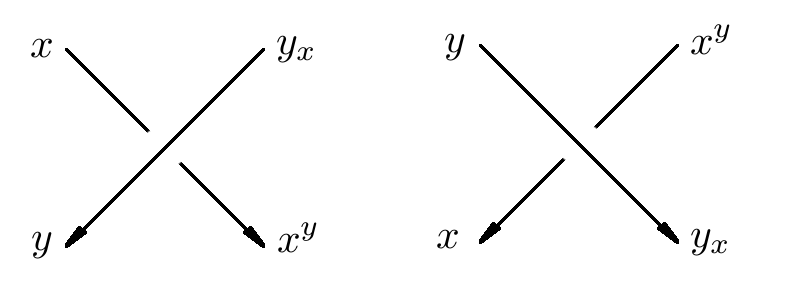}\]
(known as \textit{$X$-labelings})
correspond bijectively before and after the \textit{oriented Reidemeister 
moves.}
\[\includegraphics{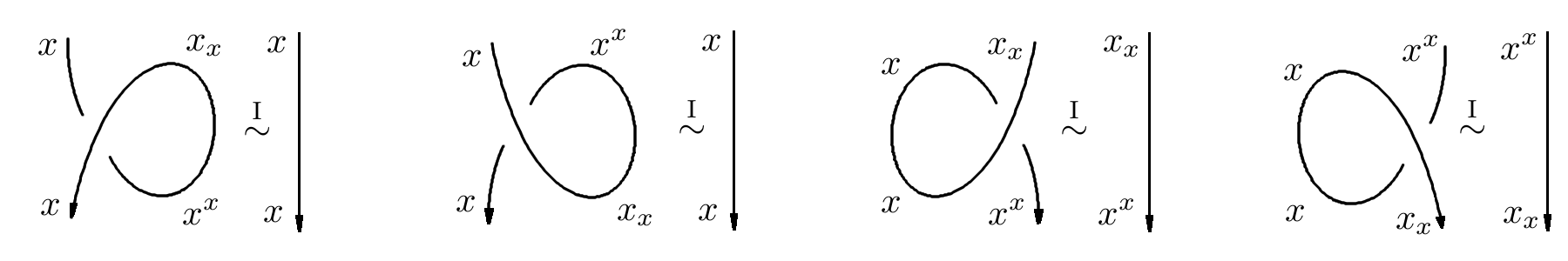}\]
\[\scalebox{0.95}{\includegraphics{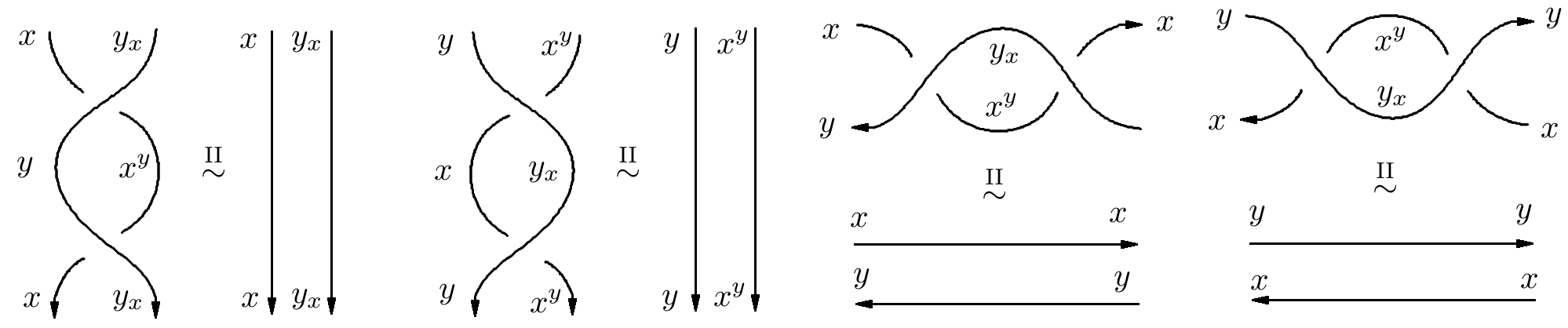}}\]
\[\includegraphics{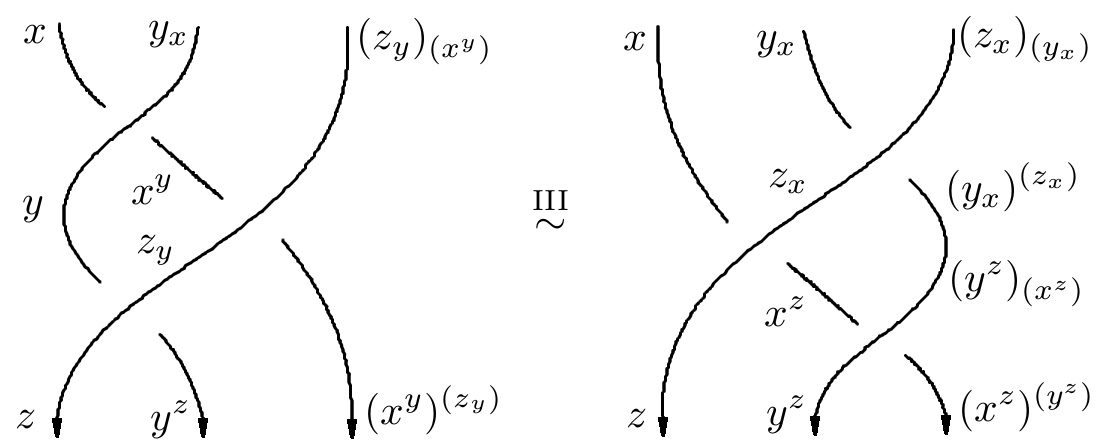}\]

In particular, by construction we have (see also \cite{KR} etc.)
\begin{theorem}
If $X$ is finite biquandle and $L$ and $L'$ are diagrams representing ambient 
isotopic oriented links, then
\[|\mathcal{L}(L,X)|=|\mathcal{L}(L',X)|\]
where $\mathcal{L}(L,X)$ is the set of $X$-labelings of $L$.
\end{theorem}
 
The number of $X$-labelings of a link diagram $L$ by a biquandle $X$ is
a computable link invariant known as the \textit{biquandle counting invariant},
denoted $\Phi_X^{\mathbb{Z}(L)}=|\mathcal{L}(L,X)|$.

\begin{example}
Consider the set $X=\{x_1,x_2\}$. There are two biquandle structures on $X$,
given by operation matrices
\[X_1=\left[\begin{array}{rr|rr}1 & 1 & 1 & 1\\2 & 2 & 2 &2 
\end{array}\right]
\quad\mathrm{and}\quad 
X_2=\left[\begin{array}{rr|rr}2 & 2 & 2 &2 \\ 1 & 1 & 1 & 1\end{array}\right]
.\]
Every classical link $L$ of $c$ components has exactly $2^c$ labelings by 
$X_1$ and $X_2$
since for any choice of base point on a component and orientation we can label 
the starting semiarc with $x_1$ or $x_2$, and then the biquandle labeling rule 
requires either keeping the label fixed at each crossing point in the 
$X_1$ case or alternating the labels in the $X_2$ case.
For example, the virtual trefoil knot $K$ below has two $X$-labelings:
\[\includegraphics{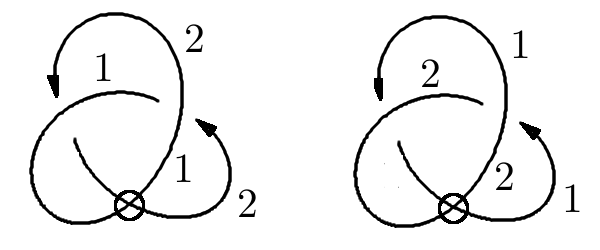}\]
For virtual links, the number of labelings by $X_2$ is $2^c$ if every component
has an even number of crossing  points (over and under) or zero if any 
component has an odd number of crossing points; for example, the 
\textit{virtual Hopf link} has no $X_2$-labelings.
\[\includegraphics{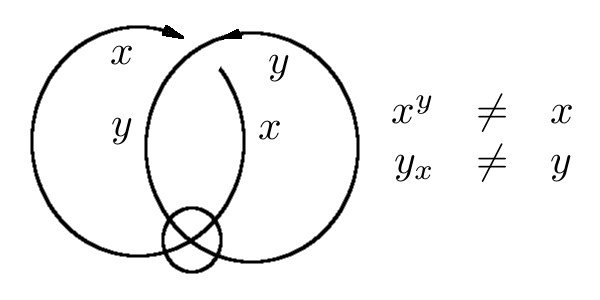}\]
As we will soon see, enhancements of the counting invariant can yield 
nontrivial information even with these simple biquandles.
\end{example}

\begin{example}
Consider the biquandle $X$ with operation matrix
\[
\left[\begin{array}{rrr|rrr}
1 & 3 & 2 & 1 & 1 & 1 \\
3 & 2 & 1 & 2 & 2 & 2 \\
2 & 1 & 3 & 3 & 3 & 3
\end{array}\right].\]
Labelings by this biquandle are the classical Fox 3-colorings of knots. Then
for example the trefoil knot $3_1$ has biquandle counting invariant
$|\mathcal{L}(3_1,X)|=9$ while the unknot $0_1$ has biquandle counting invariant
$|\mathcal{L}(0_1,X)|=3$.  
\end{example}

\section{\large\textbf{Biquandle-labeled Arrow Diagrams}}\label{A}

Let $K$ be an oriented knot diagram with a choice of base point and a 
sign and label for each crossing. The \textit{signed Gauss code} for $K$
is the ordered list of signed crossing labels encountered as we travel
around the knot from the base point following the orientation, noting whether 
we are passing over or under. For example, the oriented figure eight knot
below has signed Gauss code $U1^-O2^-U3^+O4^+U2^-O1^-U4^+O3^+$.
\[\includegraphics{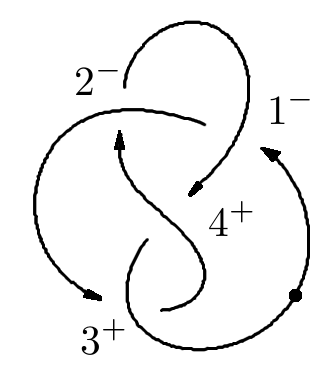}\]

Writing the Gauss code counterclockwise around a circle, we draw arrows
connecting the over-instance of each crossing to its under-instance and 
label the arrow with the crossing sign to obtain a \textit{Gauss diagram}.
\[\includegraphics{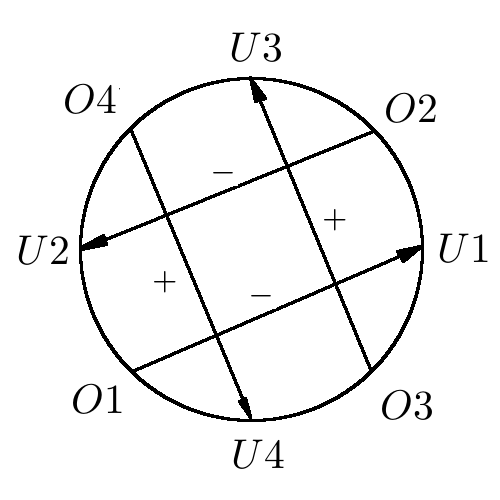}\]

In \cite{GPV}, \textit{arrow diagrams} are $\mathbb{Z}$-linear combinations
of Gauss diagrams; let us denote the set of all Gauss diagrams by $\mathcal{D}$,
so the set of all arrow diagrams is $\mathcal{A}=\mathbb{Z}[\mathcal{D}]$. A 
particular linear combination is given convenient abbreviation: a dashed arrow 
in a diagram indicates a difference of two diagrams $D-D'$ where 
\begin{itemize}
\item $D$ has an arrow as indicated by the dashed arrow,
\item $D'$ has the dashed arrow replaced with no arrow, interpreted
as a virtual crossing, and
\item $D$ and $D'$ are identical outside the neighborhood of the dashed arrow.
\end{itemize}
\[\includegraphics{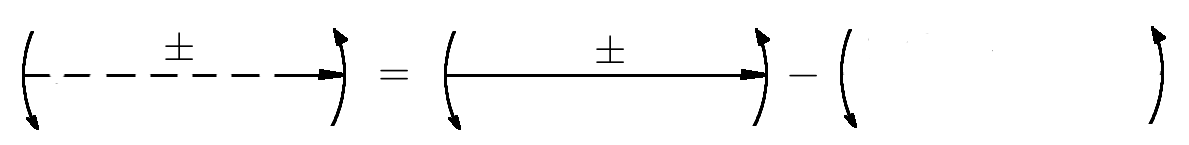}\]
Then an undashed arrow can be interpreted as a formal linear combination of 
diagrams with dashed arrows. In particular, a Gauss diagram can be expanded
as the sum of all of its subdiagrams with arrows made dashed.
\[\includegraphics{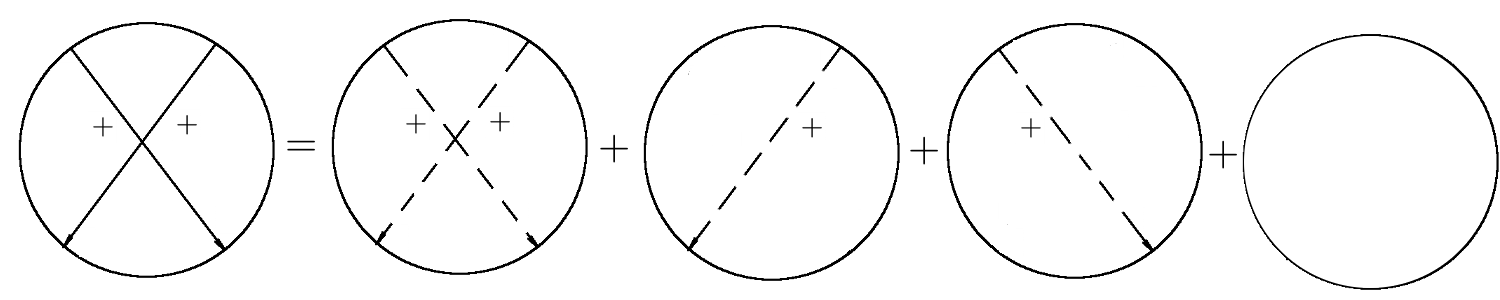}\]
We have a distinguished basis for $\mathcal{A}$
consisting of diagrams with only dashed arrows, giving us an inner product 
$\langle,\rangle:\mathcal{A}\times\mathcal{A}\to \mathbb{Z}$
by setting
\[\langle D,D'\rangle =\left\{\begin{array}{ll}
1 & D=D' \\
0 & D\ne D'
\end{array}\right.\]
for $D,D'$ basis elements and extending linearly.

The Reidemeister moves can then be interpreted as equations of elements of 
$\mathcal{A}$. The submodule $R\subset\mathcal{A}$ generated by elements of 
the form $D-D'$ where $D$ and $D'$ are Gauss diagrams related by a Reidemeister
move is called the \textit{Reidemeister submodule}. The orthogonal complement 
of this submodule in $\mathcal{A}$ is known as the \textit{Polyak Algebra}
$\mathcal{P}$. The quotient of $\mathcal{P}$ obtained by setting all diagrams 
with $\ge n$ crossings to zero is the \textit{truncated Polyak algebra of 
degree $n$}, denoted $\mathcal{P}_n$. 
Elements of $\mathcal{P}_n$ define integer-valued invariants of knots and 
links via the inner product, i.e. for $D\in\mathcal{A}$ we have
\[\phi_D(L)=\langle D,L\rangle.\]
These are known as \textit{finite type invariants of degree $n$ in the 
sense of \cite{GPV}}.

We now generalize this setup to the case of $X$-labeled Gauss diagrams.
Let $X$ be a finite biquandle and suppose $f$ is an $X$-labeling of a
knot $K$ represented by a Gauss diagram $D$. Then $f$ assigns a pair of 
elements of $X$ to each arrowhead and each arrowtail in $D$, so that every 
arrow is labeled with a sign and four elements of $X$ forming a valid 
$X$-labeling of the crossing.
\[\includegraphics{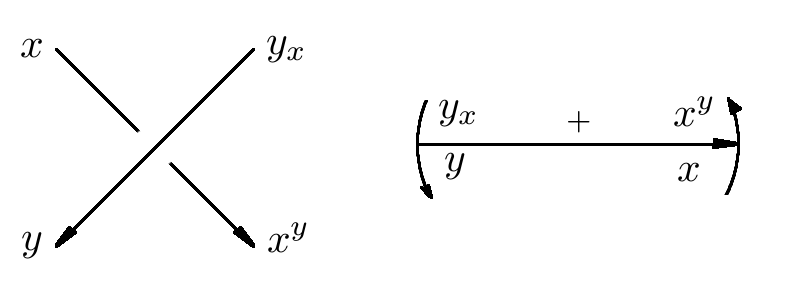}\quad
\includegraphics{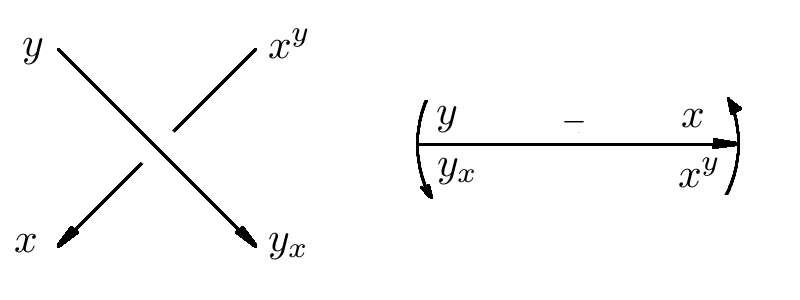}\]
Such a labeled arrow is a \textit{$X$-labeled arrow}. To keep the diagrams
as uncluttered as possible, we may list only two of the four labels 
on each arrow, one at each end, from which the other two labels can be 
recovered. 
\[\includegraphics{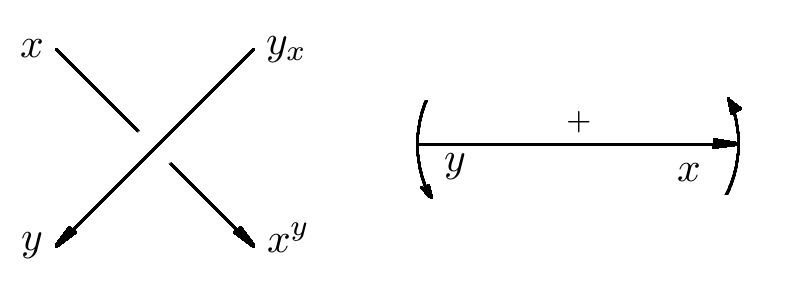}\quad
\includegraphics{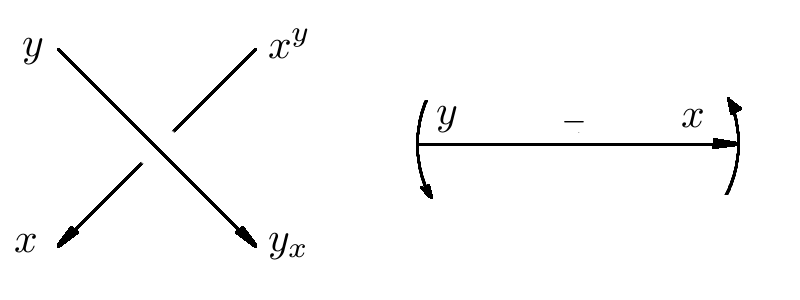}\]

\begin{definition}
A \textit{locally $X$-labeled Gauss diagram} is a Gauss diagram with pairs
of elements of $X$ at the heads and tails of each arrow determining
valid $X$-labelings of the crossings represented by the arrows. If
the two labels on the semiarcs between crossing points agree for every
semiarc, the diagram is \textit{globally $X$-labeled}.
\end{definition}

As in \cite{GPV}, we want to consider a Gauss diagram as a formal linear 
combination of subdiagrams represented as diagrams with dashed arrows. We 
observe that in such diagrams the resulting $X$-labelings are valid biquandle
labelings only locally around each arrow and in general not globally; in 
particular, the subdiagrams of a validly $X$-labeled Gauss diagram are
typically only locally $X$-labeled diagrams.

\begin{example}
Let $X$ be the biquandle with operation matrix
\[M_X=\left[\begin{array}{rr|rr}
2 & 2 & 2 & 2 \\
1 & 1 & 1 & 1
\end{array}\right].\] The $X$-labeled
virtual trefoil knot 
\[\includegraphics{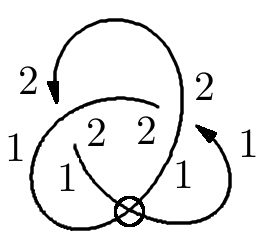}\]
has the globally $X$-labeled Gauss diagram below, which
expands as a linear combination of locally $X$-labeled arrow 
diagrams as depicted.
\[\includegraphics{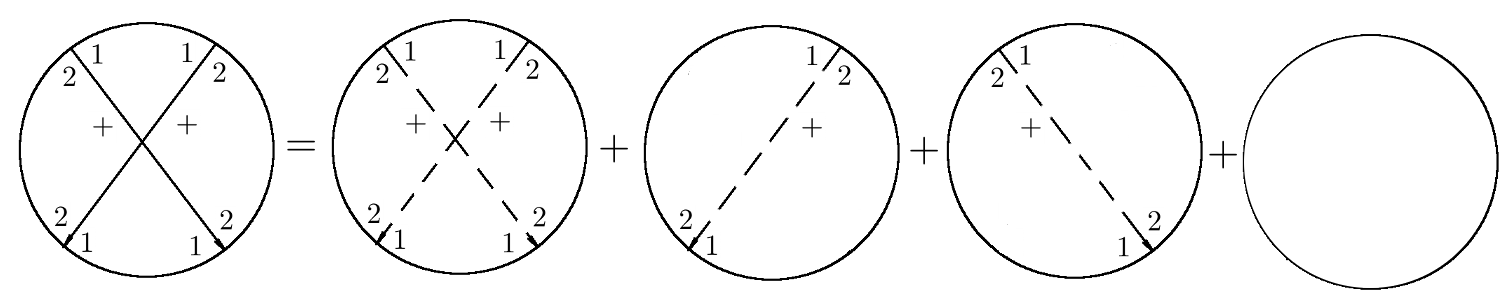}\]
\end{example}

The free abelian group generated by diagrams with locally $X$-labeled arrows
$\mathcal{A}^X$ has basis consisting of diagrams with locally $X$-labeled 
dashed arrows and inner product 
$\langle,\rangle:\mathcal{A}^X\times\mathcal{A}^X\to \mathbb{Z}$
defined as above. In particular we have:

\begin{definition}
Let $X$ be a finite biquandle. The \textit{$X$-labeled Arrow Algebra}
$\mathcal{A}^X$ is the free $\mathbb{Z}$-module generated by locally 
$X$-labeled arrow diagrams. The \textit{$X$-labeled Polyak Algebra}
$\mathcal{P}^X$ is the orthogonal complement in $\mathcal{A}^X$ of the
submodule $R^X$ generated by elements of the form $D-D'$ where $D$ and
$D'$ differ by the $X$-labeled Reidemeister moves, i.e. the relations
\begin{itemize}
\item[(i)]
\[\includegraphics{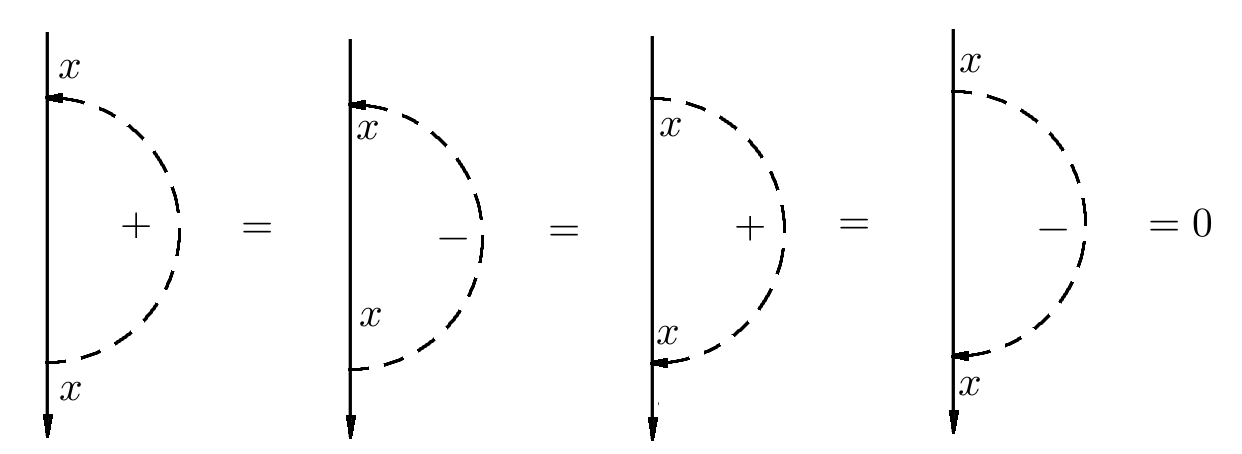}\]
\item[(ii)]
\[\includegraphics{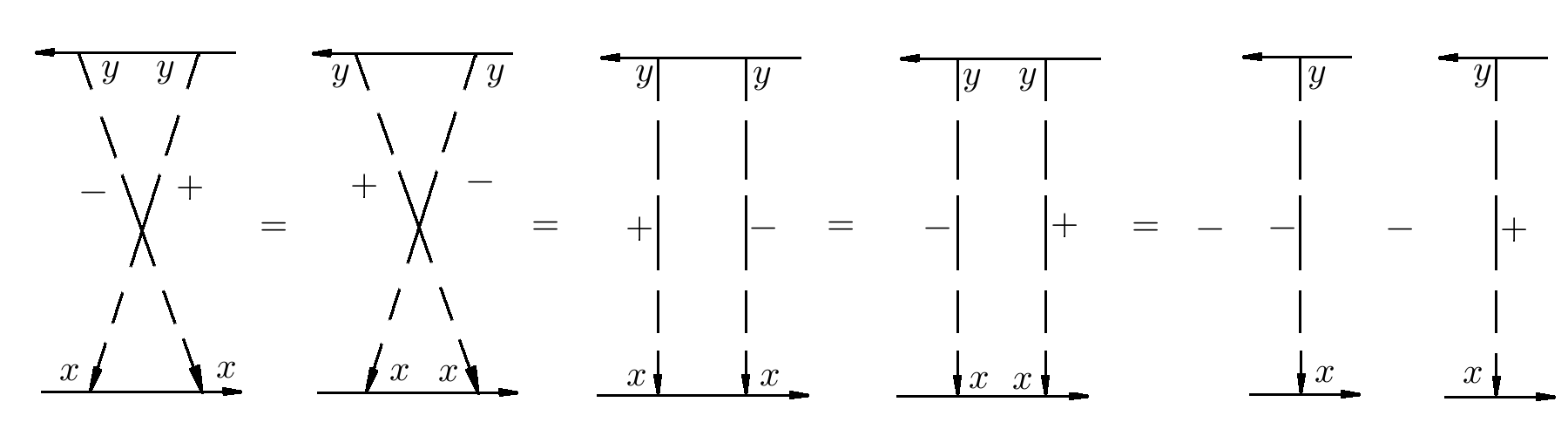}\]
\item[(iii)]
\[\includegraphics{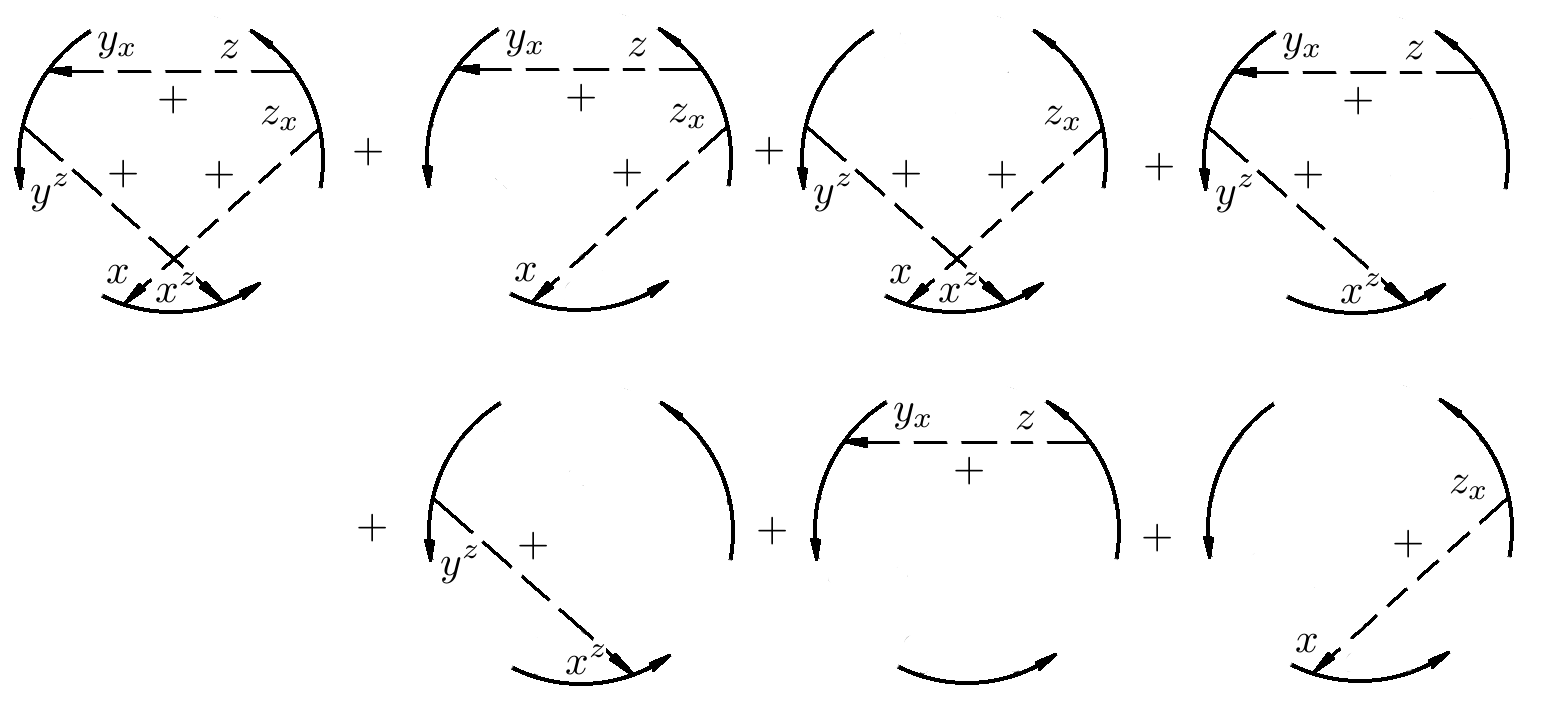}\]
\[\raisebox{1.7in}{$=$}\includegraphics{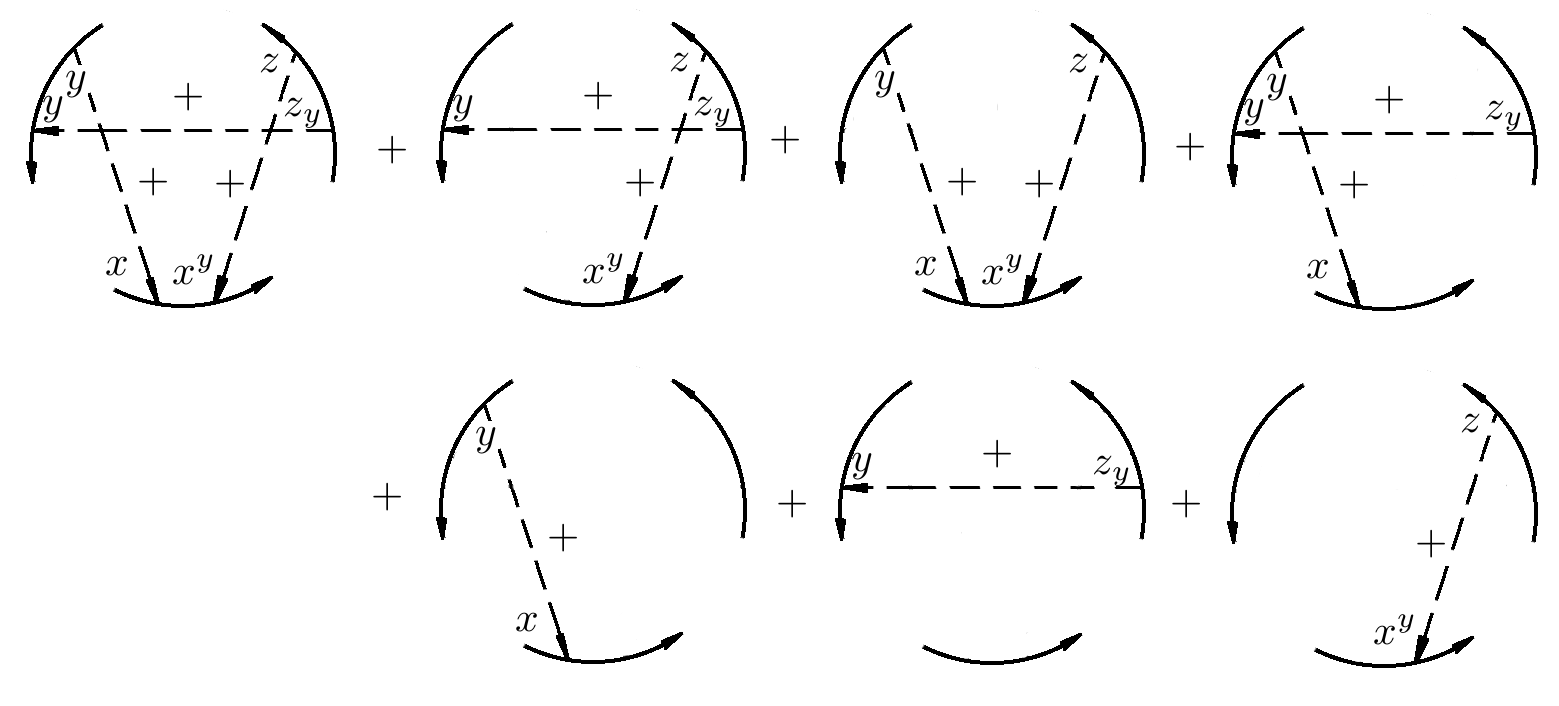}\]
\end{itemize}
for all $x,y,z\in X$. The \textit{truncated $X$-labeled Polyak Algebra} of 
degree $n$, $\mathcal{P}_n^X$, is the quotient of $\mathcal{P}^X$ obtained 
by setting all diagrams with more than $n$ arrows to zero.
\end{definition}

\begin{example}
Let $X=\{x_1\}$ be the trivial biquandle of one element. Then
the $X$-labeled arrow, Polyak and truncated Polyak algebras 
$\mathcal{A}^X$, $\mathcal{P}^X$ and $\mathcal{P}^X_n$ are isomorphic
to the original arrow, Polyak and truncated Polyak algebras defined in 
\cite{GPV}.
\end{example}

Unlike the unlabeled case, the $X$-labeled truncated Polyak algebra
with $n=1$ for a nontrivial biquandle $X$ is generally nonzero. For
any biquandle $X$, $\mathcal{A}^X_1$ is generated by $2|X|^2$ diagrams of the
form
\[\includegraphics{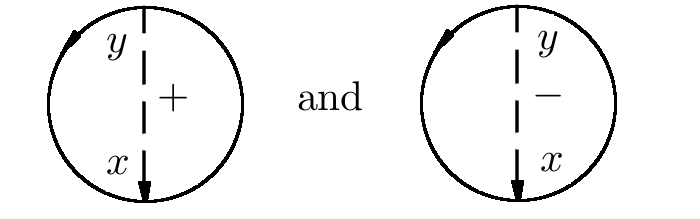}\]
for $x,y\in X$. The Reidemeister submodule is then generated by elements
of the forms
\[\includegraphics{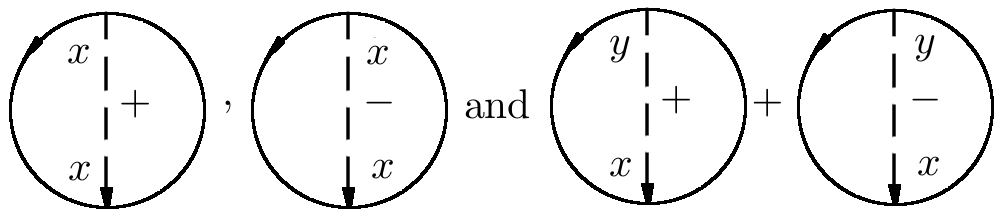}\]
and
\[\includegraphics{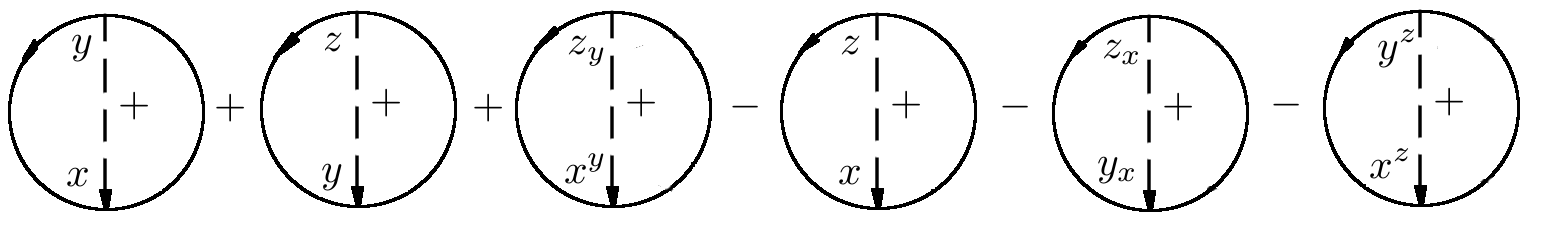}.\]
\begin{remark}
We note that the second type of generator corresponds to the biquandle 2-cocycle
condition, 
\[\phi(x,y)+\phi(y,z)+\phi(x^y,x_y)-\phi(x,z)-\phi(y_x,z_x)-\phi(x^z,y^z)=0.\]
\end{remark}

\begin{example}\label{ex:1}
Let $X$ be the biquandle with operation matrix
\[\left[\begin{array}{rr|rr}
2 & 2 & 2 & 2 \\
1 & 1 & 1 & 1
\end{array}\right].\]
$\mathcal{A}^X_1$ is generated by the eight diagrams
\[\includegraphics{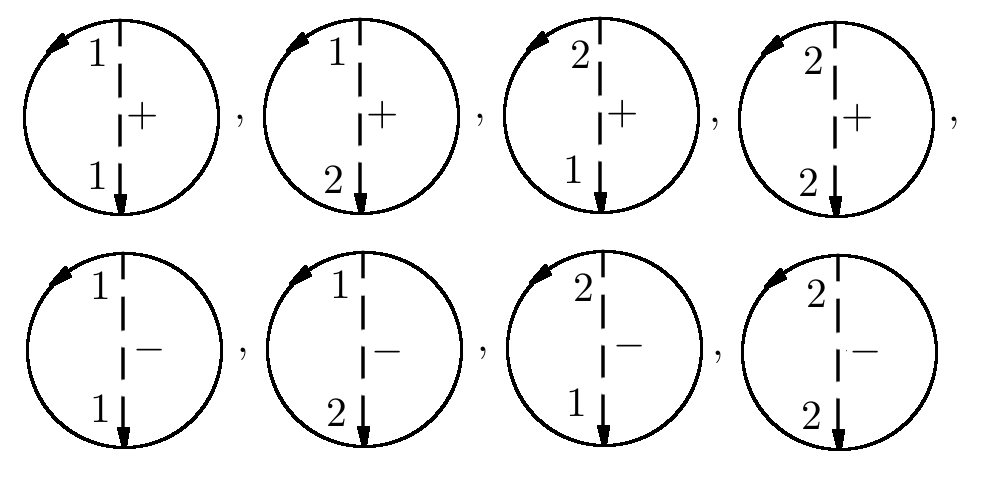}\]
with Reidemeister submodule $R^X$ generated by the six elements
\[\includegraphics{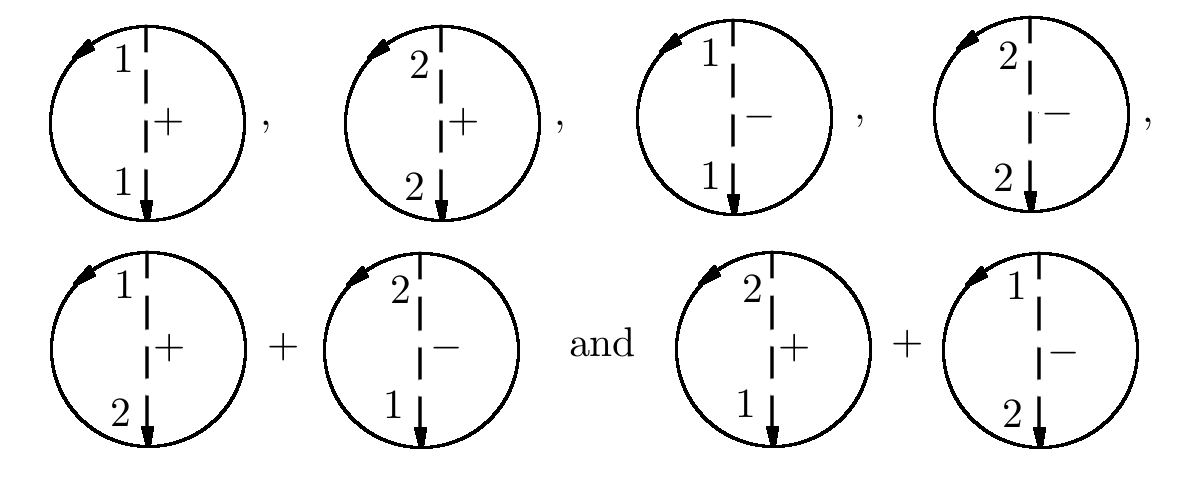}.\]
Thus, the degree 1 $X$-labeled Polyak algebra is two dimensional;
in fact it has basis
\[\left\langle \raisebox{-0.35in}{\includegraphics{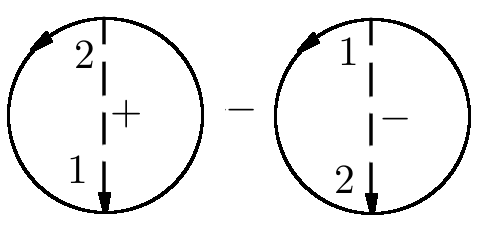}},
\raisebox{-0.35in}{\includegraphics{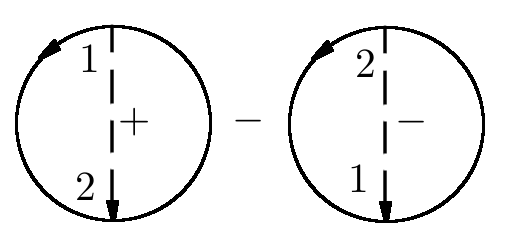}}\right\rangle.\]

\end{example}

\section{\large\textbf{Finite Type Enhancements}}\label{E}

Let $X$ be a finite biquandle and $\mathcal{P}^X_n$ the degree $n$ truncated 
$X$-labeled Polyak algebra. By construction, for any element 
$A\in\mathcal{P}^X_n$, the inner product
\[\phi_{A}(D)=\langle A,D\rangle\]
with any $X$-labeled knot diagram $D$ is unchanged by $X$-labeled Reidemeister
moves. Collecting these over a complete set of $X$-labelings of a knot
then gives us an enhancement of the biquandle counting invariant. More formally,
we have
\begin{definition}
Let $X$ be a finite biquandle and $A\in \mathcal{P}^X_n$. Then the 
\textit{finite type enhanced multiset} invariant of a knot $K$ with respect to
$X$ and $A$ is
\[\Phi_X^{A,M}(K)=\{\langle A,D\rangle \ : \ D\in\mathcal{L}(K,X)\}\]
\textit{finite type enhanced polynomial} invariant of a knot $K$ with respect 
to $X$ and $A$ is
\[\Phi_X^{A}(K)=\sum_{D\in\mathcal{L}(K,X)} u^{\langle A,D\rangle}.\]
\end{definition}

Even the simple case of the two-element biquandle in example \ref{ex:1} gives 
a nontrivial enhancement. Recall that a crossing has \textit{even parity} if 
the number of crossing points (arrow heads or tails) between its over and 
under instances is even, and \textit{odd parity} if the number of crossing 
points between its over and under crossings is odd. It is an elementary
observation to note that a crossing's parity is unchanged by Reidemeister 
moves; see \cite{VOM} for more. Moreover, a classical knot has zero crossings of
odd parity when counted algebraically, i.e., with signs equal to crossing 
signs, since the only odd parity crossings in a classical knots
are those introduced in type II moves.

\begin{example}
Let $X$ be the biquandle with two elements from example \ref{ex:1} and
\[\includegraphics{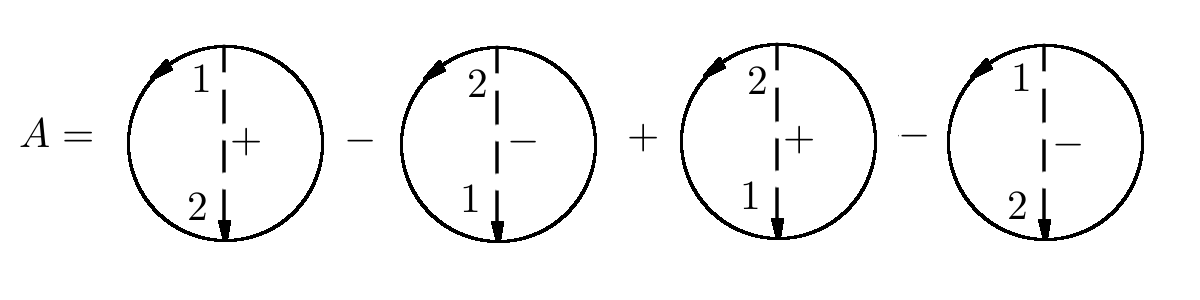}.\]
Then $\Phi_X^A(K)=2u^{P-N}$ where $P$ is the number of positively oriented 
crossings of odd parity and $N$ is the number of negatively oriented 
crossings of odd parity. For example, let $K_{3.1}$ be the virtual trefoil 
knot and let $K_U$ be the unknot,  Then even though $K_{3.1}$ and $K_U$ both 
have counting invariant value $\Phi^{\mathbb{Z}}_X=2$, we have enhancement
\[\Phi^A_X(K_{3.1})=2u^2\ne 2=\Phi_X^A(K_U).\]
Indeed, $\Phi(K)=2$ for all classical knots, so $\Phi_X^A(K)\ne 2$
implies $K$ is nonclassical.
In particular, the enhanced invariant is not determined by the counting 
invariant, and $\Phi^A_X$ is a proper enhancement.
\end{example}

\section{\large\textbf{Multicomponent Polyak Algebras}}\label{MPA}

The two element trivial biquandle, i.e. $X_1=\{x_1,x_2\}$ with operation matrix
\[\left[\begin{array}{rr|rr}
1 & 1 & 1 & 1 \\
2 & 2 & 2 & 2
\end{array}\right]\]
also has two-dimensional $\mathcal{P}^{X_1}_1$, generated by elements which look 
the same as the generators for $\mathcal{P}^{X_2}_1$ unless we list all four 
arrow labels:
\[\mathcal{P}^{X_2}_1=\left\langle \raisebox{-0.35in}{\includegraphics{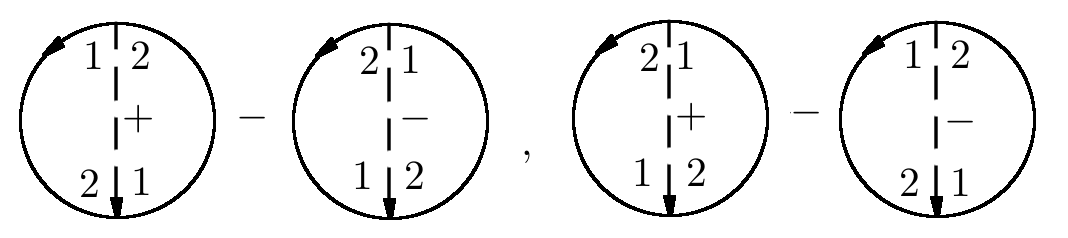}}\right\rangle\]
\[\mathcal{P}^{X_1}_1=\left\langle \raisebox{-0.35in}{\includegraphics{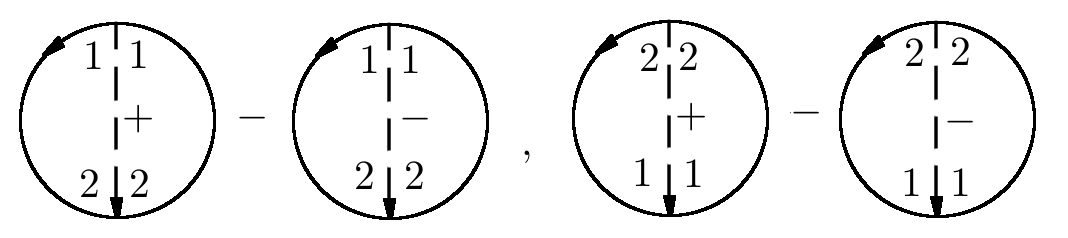}}\right\rangle.\]
Unfortunately, the enhancements defined by elements of $\mathcal{P}^{X_1}_1$
are trivial on classical and virtual knots since every $X$-labeling of a
single component virtual or classical knot is monochromatic. 

We can get something nontrivial, however, if we consider 
\textit{multi-component Polyak algebras.} A Gauss diagram of $c$ components
is a set of $c$ oriented circles with arrows connecting points on the circles
to other points on the circles. Such diagrams have a natural interpretation
as virtual link diagrams, with intra-component arrows representing 
single-component crossings and inter-component arrows representing 
multi-component crossings. For example, the 2-component Gauss diagram below
corresponds to the \textit{virtual Hopf link} depicted:
\[\includegraphics{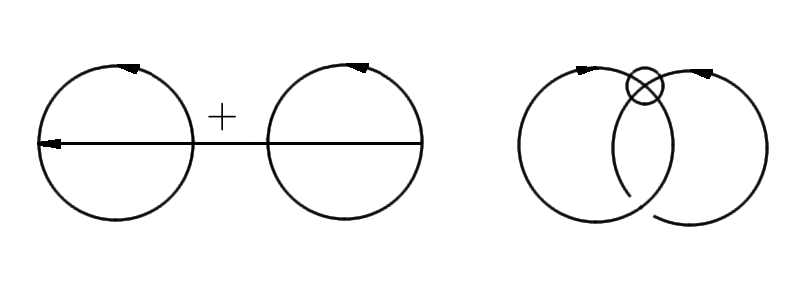}\]
We then have $c$-component Arrow, Polyak and truncated Polyak algebras
in both unlabeled and $X$-labeled varieties.

The two-component degree 1 truncated $X_1$-labeled Polyak algebra includes
the element
\[\includegraphics{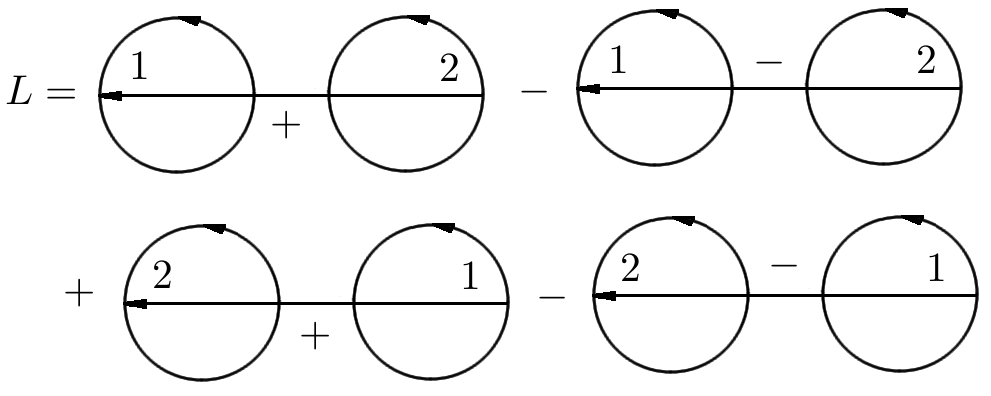}\]
For any two-component link $K$ with an even number of crossing points on each 
component, there are four $X$-labelings of $K$: two 
monochromatic and two with the components labeled with different labels.
For the monochromatic labelings we have inner product $\langle L,K\rangle=0$,
while for the non-monochronomatic labelings the inner product 
$\langle L,K\rangle$ counts the number of multicomponent crossings with
positive crossings counted with a $+1$ and negative crossings counted with a 
$-1$ -- that is, the inner product is twice the linking number of the link $K$.
Then we have:

\begin{theorem}
Let $X_1$ be the trivial biquandle of two elements, $K$ a virtual link with
an even number of crossing points on every component, and $L$ the element of 
the two-component truncated degree 1 $X_1$-labeled Polyak algebra above. Then
the finite type enhancement 
\[\Phi_x^L(K)=2+2u^{2\mathrm{lk}(K)}\]
where $\mathrm{lk}(K)$ is the linking number of $K$.
\end{theorem}

\begin{remark}
The virtual linking numbers $\mathrm{lk}_{1/2}$ and $\mathrm{lk}_{2/1}$ 
mentioned in \cite{GPV} can be recovered as finite type enhancements
by defining an ordering on the components of $K$
and extending this ordering to the circles in the diagrams generating 
$\mathcal{A}^{X_1}_1$.
\end{remark}

\section{\large\textbf{Questions}}\label{Q}

We conclude with some open questions and directions for future work, noting
that we have only scratched the surface of this subject and there remains 
much to be studied.

What is the relationship between the enhancement $\Phi^{X_2}_1$ and the arrow
sign counting invariants such as those in \cite{D,H}? What is the relationship 
between degree 1 finite type enhancements in general and biquandle cocycle 
invariants? 

We have not yet computed any $X$-labeled truncated Polyak algebras of degree
two or more; we expect the invariants defined by elements of these algebras
to be of interest.

\bigskip

\noindent
\textsc{Department of Mathematical Sciences \\
Claremont McKenna College \\
850 Columbia Ave. \\
Claremont, CA 91711} 

\end{document}